\def\IS{{\mathbb S}} 
 
\def\IK{{\mathbb K}}

\def\IC{\mathbb C} 
\def\ID{{\mathbb D}}

\def\zbar{{\overline{z}}} 
 
\def\wbar{{\overline{w}}} 
 
\def\essinf{{\rm ess\,inf}}

\documentclass[12pt]{article} 

\usepackage[psamsfonts]{amssymb} 
\usepackage{amsfonts,amsmath}  

\usepackage{epsfig,multicol}

\newtheorem{theorem}{Theorem}%[section] 
\newtheorem{lemma}{Lemma}%[section] 
\newtheorem{corollary}{Corollary}%[section] \usepackage{chngcntr}
\usepackage{apptools}
\AtAppendix{\counterwithin{lemma}{section}}

\title{Topological regularity for solutions to the generalised Hopf equation.}
\author{Gaven Martin \& Cong Yao\thanks{
Work of both authors partially supported by the New Zealand Marsden Fund.
 \newline
Institute for Advanced Study, 
Massey University,  Auckland,
New Zealand
\newline
email: g.j.martin@massey.ac.nz
}
}
\date{}
\begin{document}
\maketitle 

\begin{abstract} 
\noindent The generalised Hopf equation is the first order nonlinear equation with data $\Phi$ a holomorphic functions and $\eta\geq 1$ a positive weight,
\[ h_w\,\overline{h_\wbar}\,\eta(w) = \Phi.\]
The Hopf equation is the special case $\eta(w)=\tilde{\eta}(h(w))$ and reflects that $h$ is harmonic with respect to the conformal metric $\sqrt{\tilde{\eta}(z)}|dz|$.  This article obtains conditions on the data to ensure that a solution is open and discrete.  We also prove a strong uniqueness result.
\end{abstract}

\section{Introduction} The well known Hopf equation is the first order equation
\begin{equation}\label{h1}
h_w \overline{h_\wbar} \, \eta(h) = \Phi(w), \quad \mbox{almost every $z\in \Omega$}.  
\end{equation}
 Here $\Omega$ is a planar domain, and
 \begin{enumerate}
 \item $h:\Omega\to\IC$ lies in $W^{1,1}_{loc}(\Omega)$,
 \item $\Phi$ is holomorphic in $\Omega$,
 \item $\eta(w)\geq 1$ is locally Lipschitz (most often smooth), and
 \item the Jacobian determinant $J(w,h)\geq 0$ at almost every $w\in\Omega$.
 \end{enumerate}
For a $W^{2,1}_{loc}$ solution $h$,  differentiating (\ref{h1}) quickly reveals that $h$ satisfies the {\em tension equation},
\begin{equation}
\Delta h + 2 (\log \eta)_z(h) h_w \,h_\wbar = 0,
\end{equation}
expressing the fact that 
\[h:(\Omega,|dw|)\to  (\tilde{\Omega}, \sqrt{\eta(\tilde{w})}\, |d\tilde{w}|) \]
is a harmonic mapping, see J. Jost \cite{Jost} and G. Daskalopoulos and  R. Wentworth \cite{Wentworth} as general references for planar harmonic mappings,  or P. Duren \cite{Duren} in the case of $\eta\equiv 1$.  Elliptic regularity subsequently implies a higher degree of smoothness of $h$.  However no topological regularity can be guaranteed despite the assumption (4) above. Some aspects of the converse are well understood via Lewy's Theorem \cite{Lewy},  subsequently R. Schoen and S.T. Yau for the hyperbolic metric \cite{SY} and finally its generalisation to all metrics and degree $1$ mappings \cite{Martin}.   There  topological regularity (e.g. a homeomorphism) implies a harmonic mapping is a diffeomorphism.

\medskip

The purpose of this article is to consider the question of topological regularity in a more general setting.

\section{The generalised Hopf equation.} We consider the generalised Hopf equation:
\begin{equation}\label{1}
h_w \overline{h_\wbar} \; \eta(w)  = \Phi(w), \quad \mbox{almost every $z\in \Omega$}.  
\end{equation}
where $h$, together with the data $\Phi$ and $\eta$, satisfy (1--4) above. Of course if $h$ is harmonic to the metric $\sqrt{\tilde{\eta}}\, |d\tilde{w}|$,  then $h$ solves (\ref{1}) with $\eta(w)=\tilde{\eta}(h(w))$. Then $h$ should be locally Lipschitz to satisfy our assumptions.  This follows for the Hopf equation (\ref{h1}) by the result of  T. Iwaniec, L. V. Kovalev, J. Onninen, \cite{IKO1} who prove a much more general result which in fact shows solutions to (\ref{1}) are locally Lipschitz.

\medskip

There are natural extremal mapping problems where equations such as (\ref{1}) occur,  see \cite{Mar2} and \cite{MY} in the case $p=1$.

\medskip
We have assumed that $J(w,h)\geq 0$ to remove certain pathologies which can occur by piecing together reflections and so forth. In particular this implies $|h_\zbar|\leq |h_z|$. It follows that away from the discrete set 
\[ Z_\Phi = \{w:\Phi(w)=0\}\subset\Omega \]
we have $h_\wbar$  locally bounded and vanishing continuously on $Z_\Phi$, and hence $h\in W^{1,p}_{loc}(\Omega)$ for all $1\leq p<\infty$, and as mentioned in fact $h$ is locally Lipschitz.  This now guarantees the local integrability of the Jacobian. We have established
\begin{lemma} Let  $h:\Omega\to\IC$ be a $W^{1,1}_{loc}(\Omega)$ solution to (\ref{1}) 
with $J(w,h)> 0$ for almost every $w\in\Omega$. Then $h$ is a mapping of finite distortion.
\end{lemma} 
The theory of mappings of finite distortion can be found in \cite{AIM,HK2,IM1,IM2} which we use as basic references.  The hypothesis $J(w,h)>0$ almost everywhere is necessary. If $h(z)=a(x)$ and $\Phi=c^2>0$ is constant, $z=x+iy$, the equation becomes
\begin{equation}\label{aeqn} a'(x)\sqrt{\eta(x)}=c, \quad a(x)=a(x_0)+c\int_{x_0}^{x} \frac{dx}{\sqrt{\eta(x)}}\end{equation}
which is certainly not a mapping of finite distortion.  We will see in a moment that the case $\Phi$ constant is (at least locally) the general case.

\medskip

 A mapping $h$ is discrete if for every $\tilde{w}\in h(\Omega)$ the set $\{w\in \Omega:h(w)=\tilde{w}\}$ is a discrete subset of $\Omega$.  We will repeatedly use the following classical result whose proof we briefly sketch.
\begin{lemma} Suppose that $h:\Omega\to\IC$ is continuous in $\Omega$, open and discrete in $\Omega\setminus Z$,  where $Z$ is a discrete subset of $\Omega$. Then $h$ is open and discrete in $\Omega$.
\end{lemma}
\noindent{\bf Proof.} The Stoilow factorisation theorem implies that $h|\Omega\setminus Z =\Psi\circ f$ where $f:\Omega\setminus Z\to \Omega\setminus Z$ is a homeomorphism and $\Psi:\Omega\setminus Z\to\IC$ is holomorphic. As $h$ is continuous, $\Psi$ is locally bounded and so the points of $Z$ are removable singularities.  Hence $\Psi$ may be assumed to be holomorphic in $\Omega$. Then $f$ is continuous on $\Omega$ and one can show $f$ extends to a self homeomorphism of $\Omega$.  From this the result follows. \hfill $\Box$

\medskip The {\em Beltrami coefficient} of a mapping of finite distortion  $h$ is 
\[\mu_h = h_\zbar/h_z \]
and is defined almost everywhere. If $G$ is open and $\|\mu_h\|_{L^{\infty}(G)}=k<1$,  then $h$ is {\em quasiregular} in $G$ or {\em quasiconformal} if $h|G$ is a homeomorphism. Quasiregular mappings are open and discrete.

\section{Examples.} Here we give some examples to motivate the hypotheses we make on the data $\Phi$ and $\eta$ in equation (\ref{1}). 
\subsection{$\eta$ constant on a simply connected region $\Omega$.} Then $h_w \overline{h_\wbar}= \Phi$, a $W^{2,1}_{loc}$ solution is smooth and harmonic, and so we may find holomorphic $U,V:\Omega\to\IC$ so that $h=U+\bar V$. That $|h_\zbar|\leq |h_z|$ gives $|V'|\leq |U'|$. Away from the discrete set of points points where $\{U'=0\}$ we have $V'/U'$ holomorphic and $|V'/U'|\leq 1$, so $V'/U'$ is holomorphic in $\Omega$. The maximum principle implies that ether $|\mu_h|<1$ on $\Omega$ and $h$ is locally quasiregular - hence open and discrete - or $|\mu_h|\equiv 1$, $|U'|=|V'|$, $U'= \bar\zeta^2 V' + \bar c$, $|\zeta|=1$, and hence 
\begin{eqnarray*} h(z) & = & U+\bar V=2 \bar\zeta \Re e (\zeta U) = \bar\zeta  u(z)+c
\end{eqnarray*}
where $u$ is a real valued harmonic function.  In this case $h$ cannot be open and discrete. 

\medskip

There is no way to avoid this dichotomy without making further assumptions on $h$ such as $J(z,h)>0$ almost everywhere, or assuming that $\eta$ is nonconstant.

\subsection{$\Phi \overline{\eta_w} = |\Phi|\eta_w$.} \label{3.2}
This is a condition about the alignment of the arguments of $\Phi$ and $\eta_w^2$.   We first examine this condition in the case $\Phi\equiv 1$, so $h$ is  a solution to
\begin{equation}\label{Phi=1}
h_z \overline{h_\zbar} \; \eta(z) \equiv 1
\end{equation}  
Then the condition above is $\overline{\eta_w} = \eta_w$ which implies $\eta_y=0$, so that $\eta(z)=\eta(x)$. Then choose $a'(x)= {1}/{\sqrt{\eta(x)}} \leq 1$ and $h(z)=a(x)$ is a solution to (\ref{Phi=1}) as discussed at (\ref{aeqn}).

\medskip

 More generally we make the following calculation. Suppose that $G$ is open and we can find a holomorphic solution to
\begin{equation}\label{quad}
\phi'(z)^2 \, \Phi(\phi)\equiv 1, \quad \phi:G\to\Omega.
\end{equation} 
Then if $h_w \overline{h_\wbar}\; \hat{\eta}(w) =\Phi(w)$ we see that  $h\circ \phi$ solves (\ref{Phi=1}) with $\eta(z)=\hat{\eta}(\phi(z))$;
\[ (h\circ \phi)_z \overline{(h\circ \phi)_\zbar} \;\hat{\eta}(\phi(z)) = h_w(\phi) \overline{h_\wbar(\phi)} \; \hat{\eta}(\phi) \phi'(z)^2 = \Phi(\phi)\phi'(z)^2 =1.\]
Actually we can always solve (\ref{quad}) as soon as we can define a single valued branch of $\sqrt{\Phi}$, which we can do in a simply connected region where $\Phi\neq 0$. Then simply solve
\[   \psi'(w)=\frac{1}{\sqrt{\Phi(w)}} \neq 0 \]
 and put $\phi=\psi^{-1}$ on a maximal subdomain of definition. Next,
\[ \eta_w= \hat{\eta}_z(\phi(w)) \phi'(w) \]
Then $\eta_z  = \overline{\eta_z}$ implies 
\[ \hat{\eta}_w(\phi(z))\phi'(z)=\overline{\hat{\eta}_w(\phi(z))\phi'(z)}\]
and hence
\begin{equation}
\hat{\eta}_w(\phi(z))\Phi(\phi)=\overline{\hat{\eta}_w(\phi(z))}|\Phi(\phi)|
\end{equation}
and, writing $w=\phi(z)$  we have $\hat{\eta}_w \Phi =\overline{\hat{\eta}_w }|\Phi(w)|$  and this is the condition we are considering.

\medskip

By way of example let us consider the hyperbolic metric.  Then $\eta(w)=(1-|w|^2)^{-2}$ and $\Phi \overline{\eta_w} = |\Phi|\eta_w$ gives
$\Phi w^2 = |\Phi w^2|$, that is  $\Phi w^2 \geq 0$.  Away from this set we will find that solutions are open and discrete, see Theorem \ref{main}.

\medskip
We can rewrite the condition at \S \ref{3.2} away from the zeros of $\Phi$ as follows.
\begin{eqnarray*}
\Phi \overline{\eta_w} = |\Phi|\eta_w &\Rightarrow &
 \sqrt{\Phi} \overline{\eta_w} =  \sqrt{\bar\Phi} \eta_w  
\end{eqnarray*}
Thus $\sqrt{\Phi}  \eta_\wbar$ is real and that is if and only if $\Phi  \eta_\wbar^2 \geq 0$.
\section{Main results.} 
As we have seen,  away from $Z_\Phi$ we can arrange things locally so that $\Phi\equiv 1$.  We next examine that case.
\subsection{The special case $\Phi\equiv 1$.}
We have
\begin{equation}\label{sc1}
h_z \,\overline{h_\zbar}\, \eta(z)\equiv 1.
\end{equation}
We immediately note $0\leq \mu_h\leq 1$.   Let us remove the locally uniformly elliptic case so as  to focus on the non-classical degenerate case.
\begin{theorem} Suppose that $h$ is a $W^{1,1}_{loc}(\Omega)$ solution to (\ref{sc1}) which is locally quasiregular; $\|\mu\|_{L^\infty(G)} \leq k<1$ for any $G\Subset\Omega$.  Then on $G$ we have $h=\phi(g)$ where $g:G\to G$ is a quasiconformal diffeomorphism and $\phi:G\to\IC$ is holomorphic.
\end{theorem}
\noindent{\bf Proof.} On $G$ we have $|h_\zbar|\leq k |h_z|$. Thus $|h_z|\geq 1/\sqrt{k\eta} > 2/\sqrt{(1+k)\eta}$.  Define $A(z,\zeta)$, real valued, as follows.   
\begin{equation}
A(z,\zeta) = \left\{\begin{array}{cc} \frac{1+k}{2} |\zeta| &  |\zeta|\leq 1/\sqrt{(1+k)\eta/2} \\ 
1/(\eta|\zeta|) &  |\zeta| \geq   1/\sqrt{(1+k)\eta/2} \end{array} \right.
\end{equation}
Since $k<(1+k)/2$ and $\eta(z)$ is Lipschitz,  we may smooth $A$ along the set $|\zeta| = (1+k)/(2\eta)$ keeping the Lipschitz bound $|A(z,\zeta)-A(z,\xi)|\leq k' |\zeta-\xi| $, $k'<1$.
\medskip
\begin{center}
\scalebox{0.35}{\includegraphics{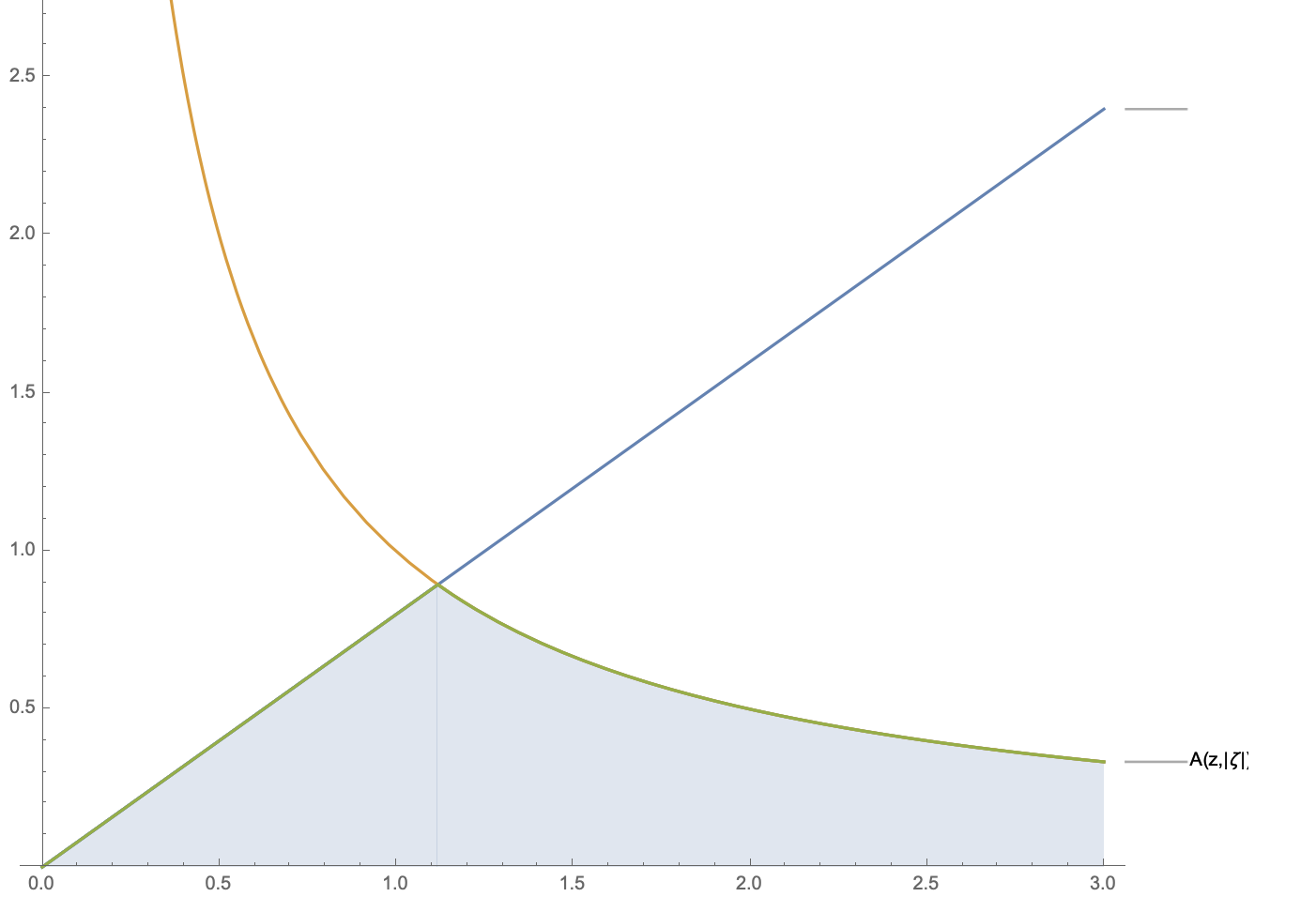}}
\end{center}
\noindent{\bf Figure 1.} The graphs $y=\frac{1+k}{2} x$ and $y=\frac{1}{\eta \, x}$. $A(z,\zeta)$ is obtained by smoothing the minimum of these two functions defined for $x\geq 0$.

\medskip
Now, simply because $|h_\zbar|  = \frac{1}{\eta(z)|h_z|}$ and so $|h_z|$ lies beyond where we have made any modification, $h$ is a solution to the following Beltrami equation with Lipschitz coefficients:
\[ h_\zbar = A(z,h_z) \frac{h_z}{|h_z|} \]
The Schauder estimates \cite{AIM}, or other more general approaches \cite{HM} show that $h\in W^{2,p}(G)$, $p<\infty$, or even $C^{\infty}(G)$ if $\eta$ is. \hfill $\Box$

\medskip

We now proceed under the additional assumption that $h\in W^{2,1}_{loc}(\Omega\setminus Z_\Phi)$ (and note that in our special case $Z_\Phi=\emptyset$).

\medskip

Differentiating (\ref{sc1}) yields
\[ 0 =h_{z\zbar} \,\mu \overline{h_z}  +  h_{z} \,\overline{h_{z\zbar}} +\frac{ \eta_\zbar }{\eta^2}=h_{z\zbar} \,\mu \overline{h_z}  + \mu h_{z} \,\overline{h_{z\zbar}} +(1-\mu) h_{z} \,\overline{h_{z\zbar}} +\frac{ \eta_\zbar }{\eta^2} \]
which simplifies to
\[ 0 = 2\mu \Re e( h_{z\zbar} \overline{h_z})+(1-\mu) h_{z\zbar} \overline{h_z} +\frac{ \eta_z}{\eta^2} \]
and the conjugate equation
\[ 0 = 2\mu \Re e( h_{z\zbar} \overline{h_z})+(1-\mu)\overline{ h_{z\zbar} }{h_z} +\frac{ \eta_\zbar}{\eta^2} \]
Subtracting the first of these equations from the second gives
\[ 0 = (1-\mu) 2i\Im m\big[h_{z\zbar} \overline{h_z}\big] +\frac{ \eta_z-\eta_\zbar}{\eta^2} \]
and hence
\begin{equation}\label{lapversion}
 \Im m\big[h_{z\zbar}\, \overline{h_z}\,\eta^2\big] =\frac{\eta_y}{ 1-\mu}
\end{equation}
While adding the same two equations gives
\[ 0 = (1+\mu) \Re e[h_{z\zbar}\, \overline{h_z}\,\eta^2] +  \eta_x  \]
This gives
\begin{equation} 
 \Im m\big[\frac{h_{z\zbar}}{h_\zbar} \big] =\frac{\eta_y}{ (1-\mu)\eta^2}, \quad   \Re e[\frac{h_{z\zbar}}{h_\zbar} ] = -\frac{\eta_x}{(1+\mu)\eta^2}
\end{equation}
This calculation has established the following lemma.
\begin{lemma}\label{deltah}
\[ h_{z\zbar} = h_\zbar \left( -\frac{\eta_x}{(1+\mu)\eta^2} + i \frac{\eta_y}{ (1-\mu)\eta^2}\right) \]
\end{lemma}
As a consequence of this formula we also have the following theorem.

\begin{theorem} Suppose we have $|\eta_y|\geq \epsilon$ and $|\eta_x|\leq M<\infty$  in $\Omega$.  Then a solution $h\in W^{2,1}_{loc}(\Omega)$ to equation (\ref{sc1})  is an open and discrete mapping.
\end{theorem}
\noindent{\bf Proof.}
As $|h_z|$ and $|h_\zbar|$ are locally bounded above and below (by a positive constant) the announced smoothness is assured by (\ref{deltah}). 
Further, we immediately see 
\[ \IK(z,h) \approx \frac{1}{1-\mu} \in L^1_{loc}(G).\]
 and so the Iwaniec-Sverak version of the Stoilow factorisation theorem \cite{IS} applies. \hfill $\Box$
 
 \medskip
 We record the following corollary related to the boundary value problem.  The proof is simply to note by topological degree theory, or a simple application of the Stoilow theorem which reduces the question to the case of mappings holomorphic in $G$ and continuous over the boundary,  that if $h:\overline{G}\to \IC$ is continuous, discrete and open in $G$, and degree $\pm1$ on the boundary $h:\partial G\to h(\partial G)$,  then it is a homeomorphism in $G$.
 
 \begin{corollary} Suppose that $\eta\geq 1$ is smooth and $\eta_y\neq 0$ in $\ID$. If $h:\overline{\ID}\to\overline{\ID}$ is a $W^{2,1}_{loc}(\ID)$ solution to (\ref{sc1}) and $h|\partial \ID\to\partial\ID$ has degree $1$,  then $h$ is a homeomorphism.
 \end{corollary}
 More generally
  \begin{corollary} Suppose that $\eta\geq 1$ is smooth in $\ID$. If $h:\overline{\ID}\to\overline{\ID}$ is a continuous $W^{2,1}_{loc}(\ID)$ solution to (\ref{sc1}), if $h|\partial \ID\to\partial\ID$ has degree $1$,  and if $h|\{\eta_y=0\}$ is injective, then $h$ is a homeomorphism.
 \end{corollary}
 \noindent{\bf Proof.} The set $\{\eta_y=0\}$ is closed and $h|\{\eta_y=0\}$ is a homeomorphism. On $\ID\setminus\{\eta_y=0\}$ $h$ is open and discrete.  The continuity of $h$ now implies the result. \hfill $\Box$
 
 \medskip
 
 As a specific consequence we note the following.
 
  \begin{corollary} Suppose that $\eta(w)=(1-|w|^2)^{-2}$ is the hyperbolic area measure in $\ID$. If $h:\overline{\ID}\to\overline{\ID}$ is a continuous $W^{2,1}_{loc}(\ID)$ solution to (\ref{sc1}), if $h|\partial \ID\to\partial\ID$ has degree $1$,  and if $h|\{y=0\}$ is injective, then $h$ is a homeomorphism.
 \end{corollary}
We need only observe that $\eta_y=2y (1-|w|^2)^{-4}$.  Actually the same will be true for any radially symmetric metric $\eta(r)$ whose derivative vanishes only at $r=0$,  since $\eta(w)=\eta(\sqrt{x^2+y^2})$ so $\eta_y= y\eta'(\sqrt{x^2+y^2})/\sqrt{x^2+y^2}$.
\subsection{The general case.}

We have seen at (\ref{quad}) how the general case we are considering,
\[ h_w\overline{h_\wbar}\,\hat{\eta} = \Phi, \]
can be reduced to the special case locally in $\Omega\setminus Z_\Phi$. Locally we have $\phi'(w)^2 \, \Phi(\phi)\equiv 1$ and 
\begin{equation}\label{geneta}
\eta(z)=\hat{\eta}(\phi(z)).
\end{equation}
Then away from $Z_\Phi$,
\begin{eqnarray*}
2i\eta_y&=&\eta_w-\eta_\wbar=2 \Im m\big(\hat{\eta}_z(\phi(z))\phi'(z)\big) \\
&=& 2 \Im m\Big(\frac{\hat{\eta}_z(\phi(z))}{\sqrt{\Phi(\phi(z))}}\Big)  = 2 \Im m\Big(\frac{\hat{\eta}_z(w)}{\sqrt{\Phi(w)}}\Big) \\
\end{eqnarray*}
Thus the problematic set is the set where $\hat{\eta}_z(w)\overline{\sqrt{\Phi(w)}}$ is real. We can express this condition more easily as $\hat{\eta}_z(w)^2 \overline{\Phi(w)} \geq 0.$

\begin{theorem}\label{main} Let $\Omega\subset \IC$  be a domain, $\Phi$ holomorphic on $\Omega$ and $Z_\Phi=\{\Phi=0\}$. Suppose that $h:\Omega\to\IC$ is  a $W^{2,1}_{loc}(\Omega\setminus Z_\Phi)$ solution to the equation
\begin{equation}\label{heqn}
h_w \overline{h_\wbar}\; \eta(w) = \Phi \quad \mbox{almost every $z\in \Omega$.}
\end{equation}
If for every $G\Subset \Omega\setminus Z_\Phi$
\[ \essinf \{\Im m(\eta_w^2\bar\Phi) : w\in G\}>0, \]
then $h$ is open and discrete.
\end{theorem}
 
\subsection{Harmonic mappings.} It is worthwhile now looking at what these conditions require in the harmonic case where $\eta(w)=\alpha(h(w))$. 
\begin{eqnarray*} \eta_w^2 \bar\Phi&=& (\alpha_z(h)h_w+\alpha_\zbar(h)\overline{h_\wbar})^2  \bar\Phi\\
 &=& \alpha_z(h)^2h_w^2 \bar\Phi+  \overline{\alpha_z(h)^2 h_\wbar^2} \bar\Phi+2 \frac{|\alpha_z(h)|^2}{\alpha(h)}|\Phi|^2 \\
\Im m(\eta_w^2 \bar\Phi)&=&\Im m\Big((\alpha_z(h)^2h_w^2  + \overline{\alpha_z(h)^2 h_\wbar^2}) \bar\Phi\Big)
 \end{eqnarray*}
 Now from (\ref{heqn}) we see $\bar\mu=|\mu|\frac{\Phi}{|\Phi|}$ and so
 \[ \Im m(\eta_w^2 \bar\Phi)=\Im m\Big(\alpha_z(h)^2h_w^2 \bar\Phi  + |\mu|^2\Phi \overline{\alpha_z(h)^2 h_w^2}\Big) \]
If we suppose that $h$ is continuously differentiable,  then $\mu$ is continuous. In a neighbourhood of a point where $|\mu|<1$ we obviously have open and discreteness.  The only interest is in the case $|\mu|=1$. However,
\begin{eqnarray*} \Im m(\eta_w^2 \bar\Phi)&=&\Im m\Big(\alpha_z(h)^2h_w^2 \bar\Phi  + \Phi \overline{\alpha_z(h)^2 h_w^2}\Big) =0
\end{eqnarray*}

\section{Uniqueness.} Here we follow the ideas of \cite{MY3}.
\begin{theorem}\label{uniq} Let $\Phi$ be holomorphic and $\eta\geq 1$ and continuous on $\ID$.  Let $g,h:\overline{\ID}\to\overline{\ID}$ be continuous mappings $g,h\in W^{1,1}_{loc}(\ID)$ with $g(0)=h(0)$ and $g(1)=h(1)$. Suppose that  $g,h \,|\IS:\IS\to\IS$ are monotone mappings of degree $1$, and that $g,h$ both solve the equation
\begin{equation}\label{99}f_w\,\overline{f_\wbar} \,\eta(w) = \Phi, \quad \mbox{almost everywhere in $\ID$} \end{equation}
If 
\[ \frac{1}{|\eta g_w h_w|- |\Phi|} \in L^{1}_{loc}(\ID\setminus Z_\Phi), \] 
then $g\equiv h$.
\end{theorem}
Notice that $|\eta(w)h_wg_w|\geq |\Phi|$ almost everywhere,  with equality holding if and only if  $J(w,h)=J(w,g)=0$.

\medskip

\noindent{\bf Proof.} The hypotheses quickly imply that $g,h$ are mappings of finite distortion in $\ID$. Let $F=g-h:\overline{\ID}\to\IC$. Then $F\in W^{1,2}_{loc}(\ID)$. Also,  as a difference of monotone maps, the total variation of the argument of $F$ on $\IS$ is at most $2\pi$, see \cite{MY3}. Away from $Z_\Phi$ we compute that
\begin{eqnarray*}
\overline{F_\zbar} & = & \overline{g_\zbar} - \overline{h_\zbar} = \frac{\Phi}{\eta}\Big(\frac{1}{g_z}-\frac{1}{h_z}\Big)\\
& = &  \frac{\Phi}{\eta g_z h_z} \, F_z= \nu\,F_z.
\end{eqnarray*}
Our hypothesis is that locally
\[ K(z,F)=\frac{1+|\nu|}{1+|\nu|}\approx \frac{1}{1-\Big|\frac{\Phi}{\eta g_z h_z}\Big|}\approx \frac{1}{|\eta g_z h_z|-|\Phi|}\]
is locally integrable away from $Z_\Phi$. Thus $F$ is a discrete open map,  \cite{IS}.  We apply the Stoilow factorisation theorem to write $F=\psi\circ f$,  where $f:\ID\to\ID$ is a homeomorphism and $\psi$ is holomorphic.  As $F$ is continuous in $\overline{\ID}$ the homeomorphism $f$ extends continuously.  As $F(0)=0$ and $F(1)=1$ the open mapping $\psi$ has degree at least two and total variation at least $4\pi$.  This contradiction establishes the result. \hfill $\Box$
 
 \begin{corollary} With the hypotheses of Theorem \ref{uniq},  if for each $z\in \ID$ either $g$ or $h$ is locally quasiregular at $z$,  then $g\equiv h$.
 \end{corollary}
 
  \begin{corollary} A locally quasiconformal solution to (\ref{99}) which is continuous in $\overline{\ID}$ is unique.
 \end{corollary}
One should compare this with Li and Tam's uniqueness theorem for quasiconformal mappings which are harmonic in the hyperbolic metric, \cite{LT}.

\end{document}